\documentclass[12pt,reqno]{amsart}
\usepackage[margin=0.85in]{geometry}
\usepackage{amsmath,amssymb,amsthm,mathtools}
\usepackage{enumitem}
\usepackage{microtype}
\usepackage[hidelinks]{hyperref}
\numberwithin{equation}{section}

\newcommand{\PSL}{\operatorname{PSL}}
\newcommand{\PGL}{\operatorname{PGL}}
\newcommand{\GL}{\operatorname{GL}}
\newcommand{\SL}{\operatorname{SL}}
\newcommand{\PSU}{\operatorname{PSU}}
\newcommand{\PSp}{\operatorname{PSp}}
\newcommand{\Aut}{\operatorname{Aut}}
\newcommand{\Out}{\operatorname{Out}}
\newcommand{\Soc}{\operatorname{Soc}}
\newcommand{\Syl}{\operatorname{Syl}}
\newcommand{\Inn}{\operatorname{Inn}}
\newcommand{\F}{\operatorname{F}}
\newcommand{\R}{\operatorname{R}}

\theoremstyle{plain}
\newtheorem{theorem}{Theorem}[section]
\newtheorem{lemma}[theorem]{Lemma}
\newtheorem{proposition}[theorem]{Proposition}
\newtheorem{corollary}[theorem]{Corollary}
\theoremstyle{remark}

\title{Products of Supersoluble Groups}
\author[J. B. Li]{Jinbao Li}
\address{School  of Mathematics and Physics, Suqian University, Suqian, Jiangsu 223800, China.}
\email{leejinbao25@hotmail.com}
\author[Y. Yang]{Yong Yang}
\address{Department of Mathematics, Texas State University, San Marcos, TX 78666, USA.}
\email{yang@txstate.edu}
\date{}
\hypersetup{pdftitle={Products of Supersoluble Groups},pdfauthor={Jinbao Li and Yong Yang},pdfsubject={Finite group factorizations},pdfkeywords={supersoluble groups, composition factors, nonsoluble length}}

\begin{document}
\begin{abstract}
Let $G=AB$ be a finite group, where $A$ and $B$ are supersoluble. We prove that every nonabelian composition factor of $G$ is isomorphic to $\PSL_2(p)$ for some prime $p\ge5$ and that the nonsoluble length of $G$ is at most one.
\end{abstract}

\maketitle
\section{Introduction}

All groups considered in this paper are finite. A classical theorem of Kegel and Wielandt states that a product of two nilpotent subgroups is soluble; see \cite[Theorem~2.4.3]{AFG}. Products of supersoluble subgroups have been studied under normality, subnormality, and permutability conditions, leading to criteria for supersolubility and nilpotence; see \cite{Baer1957,Ballester2023,GuoKondrat2015,LiLvXu2022,TangYeGuo2017} and the references therein. Without additional conditions, such a product need not be soluble. We determine the nonabelian composition factors that can occur and prove that the nonsoluble length is at most one.

\begin{theorem}\label{thm:composition}
Let $G=AB$, where $A$ and $B$ are supersoluble. Then every nonabelian composition factor of $G$ is isomorphic to $\PSL_2(p)$ for some prime $p\ge5$.
\end{theorem}

For products of two soluble subgroups, the classification in \cite{Kazarin1986} gives the family $\PSL_2(q)$ together with nine exceptional simple groups. Theorem~\ref{thm:composition} reduces this list to the groups $\PSL_2(p)$ with $p$ prime. The proof uses the reduction in \cite[Lemma~4]{KMPR2009}, the normality of a Sylow subgroup for the largest prime in a supersoluble group, and the subgroup structure of the simple groups on the list.

The \emph{nonsoluble length} $\lambda(G)$ is the minimum number of nonsoluble factors in a normal series of $G$ whose factors are either soluble or direct products of nonabelian simple groups. Bounds for this invariant involving soluble subgroups, coprime automorphisms, and subgroup factorizations were obtained in \cite{KS2015+1,KS2015+2,KS2015+3}. For products of supersoluble subgroups, we obtain the following bound.

\begin{theorem}\label{thm:main}
Let $G=AB$, where $A$ and $B$ are supersoluble. Then
\[
\lambda(G)\le1.
\]
\end{theorem}

We write $\R(G)$ for the soluble radical of $G$, its largest normal soluble subgroup, and $\Soc(G)$ for its socle, the product of its minimal normal subgroups. A group is called \emph{semisimple} if it is a direct product of nonabelian simple groups, including the trivial group as an empty product.

\begin{corollary}\label{cor:structure}
Let $G=AB$, where $A$ and $B$ are supersoluble, and let $R=\R(G)$. There exist distinct primes $p_1,\ldots,p_t\ge5$ and positive integers $r_1,\ldots,r_t$, with $t\ge0$, such that
\[
\Soc(G/R)\cong\PSL_2(p_1)^{r_1}\times\cdots\times\PSL_2(p_t)^{r_t},
\]
and $(G/R)/\Soc(G/R)$ is soluble.
\end{corollary}

\begin{proof}
Set $\overline G=G/R$. Then $\R(\overline G)=1$ and, by Theorem~\ref{thm:main}, $\lambda(\overline G)\le1$. The conclusion is immediate if $\overline G=1$. Otherwise the definition of nonsoluble length gives a semisimple normal subgroup $L$ of $\overline G$ with soluble quotient. The product of the simple direct factors in each $\overline G$-orbit is a minimal normal subgroup of $\overline G$, so $L\le\Soc(\overline G)$. Let $D$ be a minimal normal subgroup of $\overline G$. Since $\R(\overline G)=1$, the group $D$ is a direct product of nonabelian simple groups and is therefore nonsoluble. If $D\not\le L$, then minimal normality gives $D\cap L=1$, so $D\cong DL/L$ is a nonsoluble subgroup of the soluble group $\overline G/L$, a contradiction. Thus $L=\Soc(\overline G)$. Theorem~\ref{thm:composition} gives the asserted decomposition.
\end{proof}

For each prime $p\ge5$, Proposition~\ref{cf:prop:sharpness} gives a factorization of $\PGL_2(p)$ into two supersoluble subgroups with trivial intersection and one factor cyclic. This group has socle $\PSL_2(p)$ and nonsoluble length one. Thus every isomorphism type in Theorem~\ref{thm:composition} occurs, and the bound in Theorem~\ref{thm:main} is sharp. The multiplicities $r_i$ and the number of distinct primes $p_i$ in Corollary~\ref{cor:structure} are unrestricted, as the direct products following Proposition~\ref{cf:prop:sharpness} show.

Section~\ref{sec:composition} proves Theorem~\ref{thm:composition}. For the nonsoluble-length theorem, Section~\ref{sec:length-minimal} reduces a minimal counterexample to a group with socle $\PSL_2(p)^r$ and an upper simple section $\PSL_2(q)$, where $p$ and $q$ are primes. The restrictions imposed by the supersoluble factors on their actions on the socle force $p=5$ in Section~\ref{sec:prime}. The cycle bounds in Section~\ref{sec:cycles} then exclude this remaining case in Section~\ref{sec:length-a5}.

\section{Composition factors}\label{sec:composition}

For a group $X$, we write $\F(X)$ for its Fitting subgroup, $\Phi(X)$ for its Frattini subgroup, and $O_\ell(X)$ for its largest normal $\ell$-subgroup. The symbols $\Aut(X)$, $\Inn(X)$, and $\Out(X)$ denote its automorphism, inner automorphism, and outer automorphism groups. We use $\pi(X)$ for the set of prime divisors of $|X|$ and $\Syl_\ell(X)$ for the set of its Sylow $\ell$-subgroups. For a prime $\ell$ and a positive integer $n$, let $v_\ell(n)$ be the exponent of $\ell$ in $n$, let $n_\ell=\ell^{v_\ell(n)}$ be its $\ell$-part, and let $n_{2'}$ be its odd part. We denote a cyclic group of order $m$ by $C_m$ and a dihedral group of order $2m$ by $D_{2m}$. The field with $q$ elements is denoted by $\mathbb F_q$.

\subsection{Reduction to a simple direct factor}

We begin with two elementary facts about supersoluble groups.

\begin{lemma}\label{cf:lem:largest-prime}
Let $H$ be supersoluble, and let $r$ be the largest prime divisor of $|H|$. Then every Sylow $r$-subgroup of $H$ is normal in $H$.
\end{lemma}

\begin{proof}
See \cite[Chap.~VI, Theorem 9.1]{Huppert}.
\end{proof}

\begin{lemma}\label{cf:lem:affine}
Let $q=p^f$, where $p$ is a prime and $f\geq1$. Let $W$ be a nonzero additive subgroup of $\mathbb F_q$, and let $C_m\leq\mathbb F_q^{\times}$ preserve $W$ under multiplication. If $H=W\rtimes C_m$ is supersoluble, then $C_m\leq\mathbb F_p^{\times}$. In particular, $m\mid p-1$.
\end{lemma}

\begin{proof}
Since $H$ is supersoluble and $W\lhd H$ is a nontrivial $p$-group, $W$ contains an $H$-normal subgroup $L$ of order $p$. Write $L=\mathbb F_p w$ with $w\neq0$. If $\lambda\in C_m$, then $\lambda w=cw$ for some $c\in\mathbb F_p^{\times}$. Hence $\lambda=c\in\mathbb F_p^{\times}$.
\end{proof}

\begin{lemma}\label{cf:lem:component}
Suppose that
\[
G=AB=AN=BN,
\]
where $A$ and $B$ are soluble and $N$ is the unique minimal normal subgroup of $G$. Assume that $N$ is nonabelian, and write
\[
N=S_1\times\cdots\times S_t,
\qquad S_i\cong S,
\]
where $S$ is nonabelian simple. Then $A$ and $B$ are transitive on the set of components, and
\[
N\cap A=\prod_{i=1}^t(S_i\cap A),
\qquad
N\cap B=\prod_{i=1}^t(S_i\cap B).
\]
Moreover,
\begin{equation}\label{cf:eq:component-div}
|S|\mid |\Out(S)|\,|S_1\cap A|\,|S_1\cap B|.
\end{equation}
\end{lemma}

\begin{proof}
This is \cite[Lemma~4]{KMPR2009}.
\end{proof}

\begin{proposition}\label{cf:prop:minimal-reduction}
Suppose that $G=AB$ is a counterexample of minimal order to Theorem~\ref{thm:composition}. Then $G$ has a unique minimal normal subgroup
\[
N=S_1\times\cdots\times S_t,
\qquad S_i\cong S,
\]
where $S$ is nonabelian simple and is not isomorphic to $\PSL_2(p)$ for any prime $p\geq5$. Moreover, $C_G(N)=1$, $G=AN=BN$, and the supersoluble subgroups
\[
U=S_1\cap A,
\qquad V=S_1\cap B
\]
satisfy
\[
|S|\mid |\Out(S)|\,|U|\,|V|.
\]
\end{proposition}

\begin{proof}
Let $N$ be a minimal normal subgroup of $G$. The quotient
\[
G/N=(AN/N)(BN/N)
\]
is a product of two supersoluble groups and satisfies Theorem~\ref{thm:composition} by minimality. Since $G$ does not satisfy the theorem, $N$ must be nonabelian, with
\[
N=S_1\times\cdots\times S_t,
\qquad S_i\cong S,
\]
where $S$ is not isomorphic to $\PSL_2(p)$ for any prime $p\geq5$.

If $AN<G$, then
\[
AN=A(B\cap AN)
\]
is a product of two supersoluble groups and  has $\PSL_2(p)$ as a composition factor by the minimality of $G$, a contradiction. Hence $G=AN$, and similarly $G=BN$. If $M$ were a second minimal normal subgroup, then $N\cap M=1$ and the image of $N$ in $G/M$ would still have $\PSL_2(p)$ as a composition factor, again contradicting minimality. Thus $N$ is the unique minimal normal subgroup of $G$. Since $C_G(N)$ is normal and $C_G(N)\cap N=Z(N)=1$, it follows that $C_G(N)=1$. Lemma~\ref{cf:lem:component} now gives the divisibility. The subgroups $U$ and $V$ are supersoluble because they are subgroups of $A$ and $B$.
\end{proof}

Conjugation identifies $G$ with a subgroup of $\Aut(N)$. Consequently, the kernel of the action of $G/N$ on the simple direct factors of $N$ embeds in $\Out(S)^t$.

The classification in \cite{Kazarin1986} implies that the simple group $S$ in Proposition~\ref{cf:prop:minimal-reduction} is either $\PSL_2(q)$, $q>3$, or one of
\begin{equation}\label{cf:eq:kazarin-list}
\PSL_3(3),\ \PSL_3(4),\ \PSL_3(5),\ \PSL_3(7),\ \PSL_3(8),\
\PSL_4(2),\ \PSp_4(3),\ \PSU_3(8),\ M_{11}.
\end{equation}
We first consider the exceptional groups.

\subsection{The exceptional groups}

\begin{proposition}\label{cf:prop:exceptional}
No group in \eqref{cf:eq:kazarin-list} can occur as the simple group $S$ in Proposition~\ref{cf:prop:minimal-reduction}.
\end{proposition}

\begin{proof}
Write $U=S_1\cap A$ and $V=S_1\cap B$. We use \eqref{cf:eq:component-div} throughout. The normalizer structures used below are given in the \emph{Atlas} \cite{Atlas}. In the cases below, $P_r$ denotes a Sylow $r$-subgroup of $S$, chosen to contain the relevant $r$-subgroup of $U$ or $V$.

(1) Suppose first that $S=M_{11}$. We may assume that $11\mid |U|$. Since $U$ is supersoluble, its Sylow $11$-subgroup is normal, and hence
\[
U\leq N_S(P_{11})\cong 11{:}5,
\qquad |U|\mid55.
\]
As $\Out(S)=1$, equation \eqref{cf:eq:component-div} forces $144\mid |V|$. If $11\mid |V|$, then, since $V$ is supersoluble and $11$ is the largest prime divisor of $|V|$, its Sylow $11$-subgroup is normal. Hence $V\leq N_S(P_{11})\cong 11{:}5$, contradicting $144\mid |V|$. Thus $11\nmid |V|$. If $5\mid |V|$, then $5$ is now the largest prime divisor of $|V|$, so the Sylow $5$-subgroup of $V$ is normal and $V\leq N_S(P_5)$, of order $20$, again a contradiction. Hence $\pi(V)\subseteq\{2,3\}$, so a Sylow $3$-subgroup is normal in $V$. Thus
\[
V\leq N_S(P_3)\cong 3^2{:}Q_8.2,
\]
where $Q_8$ is the quaternion group of order $8$ by  \cite{Atlas}. This normalizer has order $144$, so equality holds. It is not supersoluble, since the complement acts irreducibly on its normal subgroup $3^2$, so $3^2$ is a noncyclic minimal normal subgroup. This excludes $M_{11}$.

(2) For $S=\PSL_3(3)$, a supersoluble subgroup whose order is divisible by $13$ lies in $13{:}3$. Since $|\Out(S)|=2$, the other factor, say $V$, must have order divisible by $2^3 3^2=72$. If $13\mid |V|$, then the Sylow $13$-subgroup of the supersoluble group $V$ is normal, and hence $V\leq 13{:}3$, which is incompatible with $72\mid |V|$. Thus $13\nmid |V|$. Consequently $3$ is the largest odd prime divisor of $|V|$, so a Sylow $3$-subgroup of $V$ is normal. If it has order $27$, then it is a Sylow $3$-subgroup of $S$ and its normalizer is a Borel subgroup of order $108$, whose $2$-part has order $4$. Suppose that it has order $9$, and denote it by $P$. Since the Sylow $3$-subgroups of $S$ have exponent $3$, we have $P\cong C_3^2$. Supersolubility gives a $V$-invariant subgroup of $P$ of order $3$, so $V/C_V(P)$ embeds in a line stabilizer in $\GL_2(3)$, of order $12$. We also have $|C_V(P)|_2=1$. Indeed, $S=\SL_3(3)$, and any involution of $S$ has eigenspaces of dimensions $2$ and $1$, with eigenvalues $-1$ and $1$, respectively. Its centralizer is isomorphic to $\GL_2(3)$ and has order $48$, so it cannot contain $P$. Thus
\[
|V|_2=|V/C_V(P)|_2\leq4,
\]
contradicting $8\mid|V|$.

(3) Let $S=\PSL_3(4)$. Then $|S|=2^6\cdot3^2\cdot5\cdot7$ and $|\Out(S)|=12$. After interchanging $U$ and $V$, assume that $7\mid|U|$. By Lemma~\ref{cf:lem:largest-prime}, the Sylow $7$-subgroup of $U$ is normal. It is also a Sylow $7$-subgroup of $S$, so
\[
U\leq N_S(P_7)\cong7{:}3,
\qquad |U|\mid21.
\]
Equation~\eqref{cf:eq:component-div} gives $80\mid|V|$. If $7\mid|V|$, the same normalizer argument gives $|V|\mid21$, a contradiction. Thus $5$ is the largest prime divisor of $|V|$, and its Sylow $5$-subgroup is normal. This is a Sylow $5$-subgroup of $S$, whose normalizer has order $10$ by  \cite{Atlas}. Hence $V\leq N_S(P_5)\cong D_{10}$, contradicting $80\mid|V|$.

(4) Let $S=\PSL_3(7)$. Then $|S|=2^5\cdot3^2\cdot7^3\cdot19$ and $|\Out(S)|=6$. We may assume that $19\mid|U|$. Its Sylow $19$-subgroup is normal and is a Sylow $19$-subgroup of $S$, so
\[
U\leq N_S(P_{19})\cong19{:}3,
\qquad |U|\mid57
\]
  by  \cite{Atlas}. Equation~\eqref{cf:eq:component-div} therefore gives $2^4\cdot7^3\mid|V|$. If $19\mid|V|$, the same normalizer argument gives $|V|\mid57$, a contradiction. Hence $7$ is the largest prime divisor of $|V|$, and its Sylow $7$-subgroup is normal. It has order $7^3$, so it is a Sylow $7$-subgroup of $S$. Thus $V$ lies in its Borel normalizer, whose order is
\[
7^3\frac{(7-1)^2}{(3,7-1)}=12\cdot7^3,
\]
contradicting $16\mid|V|$.

(5) Let $S=\PSL_3(8)$. Then $|S|=2^9\cdot3^2\cdot7^2\cdot73$ and $|\Out(S)|=6$. We may assume that $73\mid|U|$. By Lemma~\ref{cf:lem:largest-prime}, its Sylow $73$-subgroup is normal. It is a Sylow $73$-subgroup of $S$, so
\[
U\leq N_S(P_{73})\cong73{:}3,
\qquad |U|\mid219.
\]
Equation~\eqref{cf:eq:component-div} forces
\[
2^8\cdot7^2\mid|V|.
\]
If $73\mid|V|$, the same normalizer argument gives $|V|\mid219$, a contradiction. Thus $7$ is the largest prime divisor of $|V|$, and its Sylow $7$-subgroup $R$ is normal. It has order $7^2$, so it is a Sylow $7$-subgroup of $S$. Its normalizer is the split-torus normalizer
\[
N_S(R)\cong7^2{:}S_3,
\]
of order $294$ by  \cite{Atlas}. Hence $V\leq N_S(R)$, contradicting $2^8\mid|V|$.

(6) For $S=\PSL_4(2)\cong A_8$, choose the factor, say $U$, whose order is divisible by $7$. Then
\[
U\leq N_S(P_7),
\qquad |U|\mid21.
\]
Since $|\Out(A_8)|=2$, equation \eqref{cf:eq:component-div} forces $480\mid |V|$. If $7\mid |V|$, then $V\leq N_S(P_7)$, impossible. Hence $7\nmid|V|$, and $5$ is the largest prime divisor of $|V|$. Thus $V\leq N_S(P_5)$, whose order is $60$, again impossible.

(7) Let $S=\PSp_4(3)\cong\PSU_4(2)$. We may assume that $5\mid |U|$. Then
\[
U\leq N_S(P_5),
\qquad |N_S(P_5)|=20.
\]
Since $|S|=2^6 3^4 5$ and $|\Out(S)|=2$, equation \eqref{cf:eq:component-div} gives
\[
2^3 3^4\mid |V|.
\]
If $5\mid |V|$, then $V\leq N_S(P_5)$, a contradiction. Thus $3$ is the largest prime divisor of $|V|$, so a Sylow $3$-subgroup is normal and
\[
V\leq N_S(P_3).
\]
The latter normalizer has order $162=2\cdot3^4$, contradicting $2^3\mid|V|$. 

(8) Let $S=\PSU_3(8)$. Then $|S|=2^9\cdot3^4\cdot7\cdot19$ and $|\Out(S)|=18$. We may assume that $19\mid|U|$. Its Sylow $19$-subgroup is normal by Lemma~\ref{cf:lem:largest-prime} and is a Sylow $19$-subgroup of $S$. Hence
\[
U\leq N_S(P_{19})\cong19{:}3,
\qquad |U|\mid57.
\]
Equation~\eqref{cf:eq:component-div} forces
\[
2^8\cdot3\cdot7\mid|V|.
\]
If $19\mid|V|$, the same normalizer argument gives $|V|\mid57$, a contradiction. Thus $7$ is the largest prime divisor of $|V|$, so its Sylow $7$-subgroup $R$ is normal. The subgroup $R$ lies in a projective torus $T$ of order $21$. A generator of $R$ has a lift of order $7$ in $\operatorname{SU}_3(8)$ with the three distinct eigenvalues $\lambda,\lambda^{-1},1$, where $|\lambda|=7$, so $C_S(R)=T$. Since $R$ is characteristic in $T$, it follows that $N_S(R)=N_S(T)$. The Weyl stabilizer of this torus has order $2$, and hence $|N_S(R)|=42$ by \cite{Atlas}. Thus $V\leq N_S(R)$, contradicting $2^8\mid|V|$.

(9) It remains to consider $S=\PSL_3(5)$. Here
\[
|S|=2^5\cdot3\cdot5^3\cdot31=2(31\cdot3)(2^4\cdot5^3),
\qquad |\Out(S)|=2.
\]
We may assume that $31\mid|U|$. By Lemma~\ref{cf:lem:largest-prime}, its Sylow $31$-subgroup is normal. It is a Sylow $31$-subgroup of $S$, so
\[
U\leq N_S(P_{31})\cong31{:}3,
\qquad |U|\mid93.
\]
Equation~\eqref{cf:eq:component-div} gives $2^4\cdot5^3\mid|V|$. If $31\mid|V|$, the same normalizer argument gives $|V|\mid93$, a contradiction.  Hence $5$ is the largest prime divisor of $|V|$, and its Sylow $5$-subgroup is normal. It has order $5^3$, so it is a Sylow $5$-subgroup of $S$. Since $31\mid|U|$ and $|U|\mid93$, $|U|=31$ or $|U|=93$. If $|U|=31$, then $|V|=6000$ by Equation~\eqref{cf:eq:component-div}.  By \cite{Atlas},  $V$ involves a section of the form $S_4$ which is not supersoluble, contradicting the supersolubility of $V$. Hence we may assume $|U|=93$ and $|V|=2000$.

Let $P$ be the Sylow $5$-subgroup of $V$, which is also a Sylow $5$-subgroup of $S$. Then $P$ has the form $5_{+}^{1+2}$ with $Z(P)=[P, P]\cong C_5$ and $P/Z(P)\cong C_5^2$. By \cite{Atlas}, $V=P\rtimes  H$ with $$H=SD_{16}=\langle  x, y|x^8=y^2=1, yxy=x^3 \rangle.$$
Since $H$ acts irreducibly on $P/Z(P)\cong C_5^2$, we obtain that the chief factors of $V$ contained in $P$ is not cyclic and so $V$ is not supersoluble. This excludes the last exceptional case.

Thus the proof is complete.
\end{proof}

\subsection{The groups \texorpdfstring{$\PSL_2(q)$}{PSL(2,q)}}

Let $p$ be prime and $f$ a positive integer, and write
\[
S=\PSL_2(q),
\qquad q=p^f>3,
\qquad d=(2,q-1).
\]
For $P\in\Syl_p(S)$, the subgroup $N_S(P)$ is a Borel subgroup, that is, a point stabilizer in the natural action on the projective line $\mathbb P^1(\mathbb F_q)$. Let $T_+$ and $T_-$ be split and nonsplit tori, of orders $(q-1)/d$ and $(q+1)/d$, respectively. Their normalizers have orders
\[
|N_S(P)|=\frac{q(q-1)}d,
\qquad
|N_S(T_+)|=\frac{2(q-1)}d,
\qquad
|N_S(T_-)|=\frac{2(q+1)}d;
\]
see \cite[Chap.~II, Theorem 8.27]{Huppert}. The torus normalizers are dihedral.

\begin{lemma}\label{cf:lem:psl2-bound}
Every supersoluble subgroup of $S$ is contained in a Borel subgroup or in the normalizer of a split or nonsplit torus. Suppose that $V\leq N_S(P)$ is supersoluble, where $P\in\Syl_p(S)$, and that $V\cap P\neq1$. Let $m=|V:V\cap P|$. If $p=2$, then $m=1$. If $p$ is odd, then
\[
m\mid
\begin{cases}
p-1,& f\text{ even},\\[1mm]
(p-1)/2,& f\text{ odd}.
\end{cases}
\]
\end{lemma}

\begin{proof}
The trivial subgroup lies in every Borel subgroup. A nontrivial supersoluble subgroup $H\leq S$ has a normal subgroup $Z$ of prime order $\ell$. If $\ell=p$, then $Z$ fixes a unique point of $\mathbb P^1(\mathbb F_q)$. Since $H$ normalizes $Z$, it fixes this point and lies in its Borel stabilizer. If $\ell\ne p$, then $Z$ has exactly two fixed points on $\mathbb P^1(\mathbb F_{q^2})$. They are either both defined over $\mathbb F_q$ or interchanged by the field automorphism of $\mathbb F_{q^2}/\mathbb F_q$. The group $H$ preserves this unordered pair and therefore lies in the corresponding torus normalizer.

Under the additional hypothesis on $V$, a Hall $p'$-subgroup of $V$ has order $m$ and acts faithfully on $V\cap P$ by field multiplication. Lemma~\ref{cf:lem:affine} therefore gives $m\mid p-1$. For odd $p$, the split torus of $\PSL_2(q)$ acts on the root group through the square subgroup of $\mathbb F_q^{\times}$. Its intersection with $\mathbb F_p^{\times}$ has order $p-1$ when $f$ is even and $(p-1)/2$ when $f$ is odd. For $p=2$, $\mathbb F_2^{\times}=1$.
\end{proof}

\begin{proposition}\label{cf:prop:psl2}
If $S=\PSL_2(q)$ occurs in Proposition~\ref{cf:prop:minimal-reduction}, then $S\cong\PSL_2(p)$ for some prime $p\geq5$.
\end{proposition}

\begin{proof}
There is nothing to prove when $q$ is prime. The case $q=4$ also gives the required conclusion, since $\PSL_2(4)\cong\PSL_2(5)$. We may therefore assume that $f>1$ and $q>4$. Write $Q=G/N$, and let $K$ be the kernel of the action of $Q$ on the $t$ simple direct factors of $N$. We have $|\Out(S)|=df$.

By Lemma~\ref{cf:lem:psl2-bound}, each of $U,V$ lies in a Borel subgroup or a torus normalizer. If both lie in Borel subgroups, equation~\eqref{cf:eq:component-div} gives $(q+1)/d\mid f$. If one lies in a Borel subgroup and the other in a split-torus normalizer, it gives $(q+1)/d\mid2f$. Both are impossible because $q+1>2df$.

If both lie in torus normalizers, the defining-characteristic part of \eqref{cf:eq:component-div} gives $q\mid4f$. For odd $p$, this would imply $p^f\mid f$, which is impossible. For $p=2$, it gives $2^{f-2}\mid f$. Since $f\ge3$, this forces $f=4$: the divisibility fails for $f=3$, and $2^{f-2}>f$ for $f\ge5$. Thus $q=16$. The primes $3,5,17$ force the two torus normalizers to have opposite types, and $|S|=4080=4\cdot30\cdot34$ forces, after interchanging the factors,
\[
U\cong D_{30},\qquad V\cong D_{34}.
\]

In the natural action of $S=\PSL_2(16)=\PGL_2(16)$ on $17$ points, $U$ is the full stabilizer of an unordered pair $\{\alpha,\beta\}$. The cyclic subgroup of order $17$ in $V$ is regular. Labeling the points by $\mathbb Z/17\mathbb Z$, the reflections in $V$ have the form $z\mapsto c-z$. The reflection $z\mapsto\alpha+\beta-z$ interchanges $\alpha$ and $\beta$, so it lies in $U\cap V$. Since $(30,34)=2$, we obtain $|U\cap V|=2$ for every such pair of normalizers.

Lemma~\ref{cf:lem:component} now gives
\[
|(A\cap N)\cap(B\cap N)|=2^t,
\qquad 2^t\mid|A\cap B|.
\]
Since $G=AN=BN$, the product formula gives
\[
|A\cap B|=\left(\frac{30\cdot34}{4080}\right)^t|Q|
=\frac{|Q|}{4^t},
\]
so $8^t\mid|Q|$. But $K\leq\Out(S)^t\cong C_4^t$ and $Q/K\leq S_t$, whence
\[
v_2(|Q|)\leq2t+v_2(t!)<3t,
\]
a contradiction.

We are therefore left, after interchanging $U$ and $V$ if necessary, with
\begin{equation}\label{cf:eq:psl2-locations}
U\leq N_S(T_-),
\qquad V\leq N_S(P).
\end{equation}

We first show that $V\cap P\neq1$. If $p$ is odd and $V\cap P=1$, then $p\nmid|U|$ and \eqref{cf:eq:component-div} gives $p^f\mid2f$, which is impossible. If $p=2$ and $V\cap P=1$, then $|U|_2\leq2$ and the same divisibility gives $2^{f-1}\mid f$. This is also impossible, since $f\geq3$ and $2^{f-1}>f$.

First suppose that $p=2$. Since $|\Out(S)|=f$, Lemma~\ref{cf:lem:psl2-bound} and \eqref{cf:eq:component-div} give
\[
q(q^2-1)
\leq
f\cdot 2(q+1)\cdot q,
\]
so
\[
q-1\leq 2f.
\]
This is impossible for $f\geq3$.

Now let $p$ be odd, and write $m=|V:V\cap P|$. Equation~\eqref{cf:eq:psl2-locations} gives $|U|\leq q+1$ and $|V|\leq qm$. Since $|\Out(S)|=2f$, we have
\[
\frac{q(q-1)(q+1)}2
\leq
2f(q+1)q m,
\]
where
\[
m\leq
\begin{cases}
p-1,& f\text{ even},\\
(p-1)/2,& f\text{ odd}.
\end{cases}
\]
Hence
\begin{equation}\label{cf:eq:q-bound}
q-1\leq 4fm.
\end{equation}
If $f>1$ is odd, then
\[
1+p+\cdots+p^{f-1}\leq 2f,
\]
which is impossible. If $f$ is even and $f\geq4$, then
\[
1+p+\cdots+p^{f-1}\leq4f,
\]
which is already impossible for $p=3$. Thus $f=2$, and \eqref{cf:eq:q-bound} gives $p+1\leq8$. Therefore
\[
q\in\{9,25,49\}.
\]

If $q=25$, then
\[
|\Out(S)|\,|U|\,|V|
\mid 4\cdot26\cdot100
\]
and hence $3\nmid|\Out(S)||U||V|$, whereas $3\mid|S|$. This contradicts \eqref{cf:eq:component-div}.

Consider $q=49$. Equality is forced in all the above estimates:
\[
|U|=50,
\qquad |V|=294,
\qquad |\Out(S)|=4.
\]
The order formula gives
\[
|A\cap B|=\frac{|Q|}{4^t},
\]
so $4^t\mid |Q|$. Also,
\[
K\leq \Out(S)^t\cong (C_2\times C_2)^t.
\]
The subgroup $V$ has the form $7^2{:}6$. The component stabilizer in $B$ normalizes $V$ and hence its characteristic subgroup $O_7(V)$. Identify $O_7(V)$ with the additive group of $\mathbb F_{49}$; the complement of order $6$ is $\mathbb F_7^\times$, acting by scalars. The field involution is $\sigma:x\mapsto x^7$, so it has invariant $\mathbb F_7$-lines. By contrast, a diagonal outer automorphism acts as $x\mapsto \alpha x$ for a nonsquare $\alpha\in\mathbb F_{49}^\times$. Since every element of $\mathbb F_7^\times$ is a square in $\mathbb F_{49}^\times$, such an $\alpha$ does not lie in $\mathbb F_7$, and multiplication by $\alpha$ has no invariant $\mathbb F_7$-line. For an element in the diagonal-field coset the action is $x\mapsto \alpha x^7$ with $\alpha$ nonsquare. Its square is multiplication by
\[
N_{\mathbb F_{49}/\mathbb F_7}(\alpha)=\alpha^8,
\]
which is a nonsquare in $\mathbb F_7^\times$. Hence $x\mapsto \alpha x^7$ cannot preserve an $\mathbb F_7$-line, since on such a line its square would act as a square scalar from $\mathbb F_7^\times$. The component stabilizer in $B$ is supersoluble, so its outer image in each coordinate is contained in the subgroup generated by the field involution. Since $G=BN$, elements of $K$ lift to elements of $B$ normalizing every component. Hence
\[
|K|_2\leq 2^t.
\]
Since $Q/K\leq S_t$,
\[
v_2(|Q|)\leq t+v_2(t!)<2t,
\]
contrary to $4^t\mid |Q|$.

Finally let $q=9$, so $S\cong A_6$. The divisibility \eqref{cf:eq:component-div} and the supersoluble subgroup bounds give, after interchanging the factors if needed, the three possibilities
\[
(|U|,|V|)=(5,18),\quad (10,9),\quad (10,18).
\]
For the first two pairs, $|U||V|=90$, and hence
\[
4^t\mid |Q|.
\]
Let $P$ be a Sylow $3$-subgroup of $A_6$, so $P\cong C_3^2$. Its normalizer in $\Aut(A_6)$ has shape $3^2{:}SD_{16}$, where $SD_{16}$ is the semidihedral group of order $16$, whereas $N_{A_6}(P)/P\cong C_4$. This $C_4$ is the Frattini subgroup $\Phi(SD_{16})$, the intersection of the maximal subgroups of $SD_{16}$, and $SD_{16}/\Phi(SD_{16})\cong C_2^2$. Hence, if $H\leq SD_{16}$ maps onto the full quotient, then $H\Phi(SD_{16})=SD_{16}$; by the defining property of the Frattini subgroup, this forces $H=SD_{16}$. Thus a subgroup of the outer normalizer which maps onto the full group $\Out(A_6)\cong C_2^2$ contains the whole $SD_{16}$-action on $P$, and that action is irreducible. The normalizer structure is the one given for $A_6$ and its automorphism groups in the \emph{Atlas} \cite{Atlas}. A supersoluble component stabilizer containing $P$ must preserve a subgroup of order $3$, so it cannot induce the full outer group. Therefore the kernel $K$ of the component action satisfies $|K|_2\leq2^t$, and again
\[
v_2(|Q|)\leq t+v_2(t!)<2t,
\]
a contradiction.

For the remaining pair $(10,18)$, one has locally
\[
U\cong D_{10},
\qquad V\cong 3^2{:}2.
\]
We determine the intersection size for these two conjugacy classes. There are $36$ Sylow $5$-subgroups in $A_6$, and hence $36$ subgroups $D_{10}=N_{A_6}(P_5)$. There are $10$ Sylow $3$-subgroups $P\cong C_3^2$; for each such $P$, the quotient $N_{A_6}(P)/P\cong C_4$ has a unique subgroup of order $2$, so there are $10$ subgroups of shape $3^2{:}2$ of the required type. The involutions of $A_6$ are the $45$ double transpositions. Each $D_{10}$ contains $5$ involutions, so each involution lies in
\[
\frac{36\cdot5}{45}=4
\]
such subgroups. Each subgroup $3^2{:}2$ contains $9$ involutions, so each involution lies in
\[
\frac{10\cdot9}{45}=2
\]
such subgroups. Consequently the number of triples $(D,V,z)$, where $D\cong D_{10}$, $V\cong3^2{:}2$, and $z$ is an involution in $D\cap V$, is
\[
45\cdot4\cdot2=360.
\]
On the other hand, there are exactly $36\cdot10=360$ ordered pairs $(D,V)$. Such a pair contains at most one common involution: the product of two distinct involutions in $D_{10}$ has order $5$, whereas the product of two distinct involutions in $3^2{:}2$ lies in the normal $3^2$ and has order $3$. Thus every ordered pair has exactly one common involution. In particular,
\[
|U\cap V|=2.
\]
Hence
\[
2^t\mid |(A\cap N)\cap(B\cap N)|\mid |A\cap B|.
\]
Since now
\[
|A\cap B|=\frac{|Q|}{2^t},
\]
we again obtain $4^t\mid |Q|$, and the same outer-kernel argument gives the same contradiction.

Thus the only possible nonprime field size is $q=4$, which gives $S\cong\PSL_2(5)$.
\end{proof}

\begin{proof}[Proof of Theorem~\ref{thm:composition}]
Suppose that $G=AB$ is a counterexample of minimal order, and let $S$ be the simple group given by Proposition~\ref{cf:prop:minimal-reduction}. By \cite{Kazarin1986}, either $S\cong\PSL_2(q)$ or $S$ belongs to \eqref{cf:eq:kazarin-list}. Proposition~\ref{cf:prop:exceptional} excludes the latter possibility. Proposition~\ref{cf:prop:psl2} then gives $S\cong\PSL_2(p)$ for a prime $p\geq5$, contrary to the choice of $S$.
\end{proof}

Every simple group in Theorem~\ref{thm:composition} occurs as the socle of a product of two supersoluble subgroups, even when the factors have trivial intersection and one of them is cyclic.

\begin{proposition}\label{cf:prop:sharpness}
For every prime $p\geq5$, the group $G=\PGL_2(p)$ has a factorization $G=AB$ with
\[
A\cong C_p{:}C_{p-1},\qquad B\cong C_{p+1},\qquad A\cap B=1.
\]
In particular, $A$ and $B$ are supersoluble, $\Soc(G)=\PSL_2(p)$, and $\lambda(G)=1$.
\end{proposition}

\begin{proof}
Let $A$ be a point stabilizer in the natural action of $G$ on the projective line. Then $A\cong C_p{:}C_{p-1}$. Identify the underlying two-dimensional vector space over $\mathbb F_p$ with $\mathbb F_{p^2}$. Multiplication by $\mathbb F_{p^2}^{\times}$ induces a cyclic subgroup $B$ of $G$ of order $p+1$, acting regularly on the projective line. Hence $G=AB$ and $A\cap B=1$. Both factors are supersoluble, and $\Soc(G)=\PSL_2(p)$. The nonabelian simple normal subgroup $\PSL_2(p)$ has index $2$ in $G$, so $\lambda(G)=1$.
\end{proof}

For distinct primes $p_1,\ldots,p_t\geq5$ and positive integers $r_1,\ldots,r_t$, choose factorizations $G_i=A_iB_i$ with $G_i=\PGL_2(p_i)$ as above. Then
\[
\prod_{i=1}^tG_i^{r_i}
=
\left(\prod_{i=1}^tA_i^{r_i}\right)
\left(\prod_{i=1}^tB_i^{r_i}\right)
\]
is a product of two supersoluble groups with socle
\[
\prod_{i=1}^t\PSL_2(p_i)^{r_i}.
\]
Thus neither the number of simple direct factors nor the number of their isomorphism types is bounded in Corollary~\ref{cor:structure}.

\section{Structural reductions}\label{sec:structural}\label{sec:length-minimal}

\begin{lemma}\label{lem:local-psl2}
Let $S=\PSL_2(p)$, where $p\ge5$ is prime, and let $U,V\le S$ be supersoluble subgroups such that
\[
|S|\mid 2|U|\,|V|.
\]
After interchanging $U$ and $V$, we have
\[
p\mid |U|,\qquad p\nmid |V|,\qquad \frac{p+1}{2}\mid |V|.
\]
With this choice, the following hold.
\begin{enumerate}
\item[(i)] At least one of $U,V$ is nonnilpotent.
\item[(ii)] There exists $P\in\Syl_p(S)$ such that
\[
U\le N_S(P)\cong C_p:C_{(p-1)/2},
\qquad
(|V|,p-1)_{2'}=1.
\]
\item[(iii)] The subgroup $V$ is contained in the normalizer of a nonsplit torus of $S$. In particular, $|V|\mid p+1$.
\end{enumerate}
\end{lemma}

\begin{proof}
The divisibility hypothesis implies that $p$ divides at least one of
$|U|,|V|$. After interchanging the subgroups, assume that $p\mid|U|$.
Since $p$ is the largest prime divisor of $|S|$, the Sylow $p$-subgroup
$P$ of $U$ is normal. Hence
\[
U\le N_S(P)\cong C_p:C_{(p-1)/2}.
\]
If $p\mid|V|$ as well, the same argument gives $|V|\mid p(p-1)/2$.
Then
\[
|S|\mid2\left(\frac{p(p-1)}2\right)^2,
\]
which would imply $p+1\mid p(p-1)$, impossible for $p\ge5$.
Thus $p\nmid|V|$. Since $|U|\mid p(p-1)/2$, the divisibility hypothesis
also gives
\begin{equation}\label{eq:length-local-2}
p+1\mid2|V|.
\end{equation}

By Lemma~\ref{cf:lem:psl2-bound}, $V$ lies in a Borel subgroup or a
split- or nonsplit-torus normalizer. In the Borel case, $p\nmid|V|$
gives $|V|\mid(p-1)/2$, contradicting \eqref{eq:length-local-2}.
In the split-torus case, $|V|\mid p-1$, so
\eqref{eq:length-local-2} would give $p+1\mid2(p-1)$, also impossible.
Hence $V\le N_S(T)$ for a nonsplit torus $T$, and $|V|\mid p+1$.
Since $(p-1,p+1)=2$, this also gives $(|V|,p-1)_{2'}=1$.

If both $U$ and $V$ are nilpotent, then $C_S(P)=P$ forces $U=P$.
Consequently
\begin{equation}\label{eq:length-local-1}
p^2-1\mid4|V|.
\end{equation}
Since $|V|\le p+1$, we obtain $p-1\le4$, so $p=5$.
Then $N_S(T)\cong S_3$, while \eqref{eq:length-local-1} gives
$6\mid|V|$. Thus $V=N_S(T)\cong S_3$, contradicting nilpotence.
\end{proof}

For $N\lhd G=AB$, the subgroup $AN\cap BN$ is called the factorizer of $N$ in this factorization. The next lemma determines its action on the simple factors of $N$.

\begin{lemma}\label{lem:factorizer}
Let $G=AB$, with $A$ and $B$ supersoluble, and let
\[
N=S_1\times\cdots\times S_r\lhd G,
\]
where the $S_i$ are nonabelian simple groups. Let
\[
F=AN\cap BN,\qquad A_0=A\cap BN,\qquad B_0=B\cap AN.
\]
Then $F$ normalizes every simple direct factor $S_i$ of $N$.
\end{lemma}

\begin{proof}
By Dedekind's identity,
\[
F=A_0N=B_0N.
\]
Moreover, if $f=ab\in F$, with $a\in A$ and $b\in B$, then $a=fb^{-1}\in BN$ and $b=a^{-1}f\in AN$. Hence $a\in A_0$, $b\in B_0$, and
\[
F=A_0B_0=A_0N=B_0N.
\]

Let $\Delta$ be an $F$-orbit on the set of simple factors of $N$. Set
\[
M=\prod_{S_i\in\Delta}S_i,
\qquad
L=\prod_{S_i\notin\Delta}S_i.
\]
Then $M,L\lhd F$ and $N=M\times L$. In $F/L$, let $\widetilde M$ be the image of $M$. The group $F/L$ is transitive on the simple factors of $\widetilde M$, so $\widetilde M$ is a minimal normal subgroup of $F/L$.

Now let
\[
C=C_{F/L}(\widetilde M)
\]
and let $H=(F/L)/C$. Since $Z(\widetilde M)=1$, we have $C\cap\widetilde M=1$, so $\widetilde M$ embeds in $H$. Moreover
\[
C_H(\widetilde M)=1.
\]
Hence $\widetilde M$ is the unique minimal normal subgroup of $H$: a second minimal normal subgroup would intersect $\widetilde M$ trivially and therefore centralize it.

Let $\widehat A$ and $\widehat B$ be the images of $A_0$ and $B_0$ in $H$. Then
\[
H=\widehat A\widehat B
 =\widehat A\widetilde M
 =\widehat B\widetilde M.
\]
The groups $\widehat A$ and $\widehat B$ are supersoluble. Thus Lemma~\ref{cf:lem:component} applies. Writing
\[
\widetilde M=T_1\times\cdots\times T_m,
\]
we obtain transitivity of both $\widehat A$ and $\widehat B$ on the $T_i$, together with
\[
\widetilde M\cap\widehat A
 =\prod_{i=1}^m(T_i\cap\widehat A),
\qquad
\widetilde M\cap\widehat B
 =\prod_{i=1}^m(T_i\cap\widehat B),
\]
and
\[
|T_1|\mid |\Out(T_1)|\,|T_1\cap\widehat A|\,|T_1\cap\widehat B|.
\]
By Theorem~\ref{thm:composition}, $T_1\cong\PSL_2(p)$ for a prime $p\ge5$. Lemma~\ref{lem:local-psl2}(i) shows that at least one of
\[
U=T_1\cap\widehat A,
\qquad
V=T_1\cap\widehat B
\]
is nonnilpotent. Suppose that $U$ is nonnilpotent.

Let
\[
Y=\widetilde M\cap\widehat A
  =U_1\times\cdots\times U_m,
\]
where $U_i=T_i\cap\widehat A$ and all $U_i$ are isomorphic to $U$. Since $Y\lhd\widehat A$ and a finite supersoluble group $X$ satisfies $X'\le \F(X)$,
\[
[Y,\widehat A]\le Y\cap\widehat A'
   \le Y\cap\F(\widehat A)
   \le \F(Y).
\]
Hence $\widehat A$ centralizes
\[
Y/\F(Y)\cong\prod_{i=1}^m U_i/\F(U_i).
\]
Each factor is nontrivial because $U_i$ is nonnilpotent. A nontrivial permutation of the $U_i$ acts nontrivially on this direct product. Therefore $\widehat A$ normalizes every $U_i$. Since $\widehat A$ is transitive on them, $m=1$.

Thus every $F$-orbit $\Delta$ has size one, as required.
\end{proof}

Assume that Theorem~\ref{thm:main} is false and choose a counterexample $G=AB$ of least order. Since
\[
\lambda(G)=\lambda(G/\R(G)),
\]
the quotient by a nontrivial soluble radical would be a smaller counterexample. Hence
\[
\R(G)=1.
\]
Let
\[
E=\Soc(G).
\]
Then $E$ is a direct product of nonabelian simple groups and
\[
C_G(E)=1.
\]
Indeed, if $C_G(E)\ne1$, then this normal subgroup of $G$ contains a minimal normal subgroup of $G$, which lies in $E\cap C_G(E)=Z(E)=1$.

Let $\Omega_E$ be the set of simple direct factors of $E$, and let $K$ be the kernel of the action of $G$ on $\Omega_E$. Since each component is $\PSL_2(p)$ with $p$ prime,
\[
K/E\hookrightarrow\prod_{S\in\Omega_E}\Out(S)
\]
is an elementary abelian $2$-group. If $G/K$ were soluble, then $G$ would have nonsoluble length at most one. Hence $G/K$ is nonsoluble.

Let
\[
L/K=\R(G/K).
\]
Choose a minimal normal subgroup $M/L$ of the nontrivial group $G/L$.
Since $\R(G/L)=1$, this subgroup is nonabelian, so
\[
M/L\cong T^t,
\qquad
T\cong\PSL_2(q)
\]
for a prime $q\ge5$. Set
\[
H=M/K,\qquad R=L/K.
\]
Then $R$ is the soluble radical of $H$, and $H/R\cong T^t$.

\begin{lemma}\label{lem:M-two-layers}
With the preceding notation, $\lambda(M)\ge2$.
\end{lemma}

\begin{proof}
Since $M\lhd G$ and $\R(G)=1$, we have $\R(M)=1$. We claim that $\Soc(M)\le L$. Let $C$ be a minimal normal subgroup of $M$. If $C\not\le L$, then $C\cap L=1$. Since $E\le K\le L$ and $E\lhd M$,
\[
[C,E]\le C\cap E\le C\cap L=1.
\]
Thus $C\le C_G(E)=1$, a contradiction. Hence every minimal normal subgroup of $M$ lies in $L$.

If $\lambda(M)\le1$, then, because $\R(M)=1$, the quotient $M/\Soc(M)$ is soluble. Since $\Soc(M)\le L$, this would make $M/L$ soluble, contrary to the choice of $M$.
\end{proof}

\begin{proposition}\label{prop:upper-factorized}
With the preceding notation, we have
\[
G=AM=BM.
\]
Moreover $t=1$, and
\[
\PSL_2(q)\le G/L\le\PGL_2(q).
\]
\end{proposition}

\begin{proof}
Let
\[
X=AM\cap BM.
\]
The factorizer identity gives
\[
X=(A\cap BM)(B\cap AM),
\]
so $X$ is again a product of two supersoluble subgroups and $M\lhd X$. If $X<G$, minimality of $G$ gives $\lambda(X)\le1$. Intersecting a suitable normal series of $X$ with $M\lhd X$ gives $\lambda(M)\le\lambda(X)$, contrary to Lemma~\ref{lem:M-two-layers}. Hence $X=G$, and therefore $G=AM=BM$.

Apply Lemma~\ref{lem:factorizer} in $G/L$ to the semisimple normal subgroup $M/L=T^t$. Its factorizer is all of $G/L$, so $G/L$ fixes every simple factor of $T^t$. Minimal normality of $T^t$ requires transitivity on these factors, and hence $t=1$.

Suppose that $D/L$ is a second minimal normal subgroup of $G/L$.
Since $\R(G/L)=1$, the group $D/L$ is nonabelian semisimple and hence
nonsoluble. We have
\[
(D/L)\cap(M/L)=1,
\]
so $D\cap M=L$ and
\[
DM/M\cong D/L.
\]
Thus the supersoluble group
\[
G/M\cong A/(A\cap M)
\]
has a nonsoluble normal subgroup, a contradiction. Hence $T=M/L$ is the unique minimal normal subgroup of $G/L$. Since $\R(G/L)=1$, its centralizer is trivial, so
\[
T\le G/L\le\Aut(T)=\PGL_2(q).\qedhere
\]
\end{proof}

\begin{lemma}\label{lem:active-lower}
The group $G$ has a unique minimal normal subgroup
\[
N=E=\Soc(G)\cong S^r,
\qquad S\cong\PSL_2(p),
\]
where $p\ge5$ is prime. The group $H=M/K$ acts faithfully on the simple
factors of $N$, and its permutation image is nonsoluble. Moreover,
\[
C_G(N)=1,
\qquad K/N\le\Out(S)^r\cong C_2^r.
\]
\end{lemma}

\begin{proof}
Write
\[
E=N_1\times\cdots\times N_u
\]
as the product of the minimal normal subgroups of $G$. Since $K$ is the kernel of the action on all simple components of $E$ and $K\le M$, the group $H=M/K$ acts faithfully on the union of these component sets. If the image of $H$ on the components of every $N_i$ were soluble, then $H$ would embed in a direct product of soluble groups and would itself be soluble. This is impossible because
\[
H/R\cong\PSL_2(q).
\]
Choose $N=N_i$ such that the group induced by $H$ on its simple factors
is nonsoluble. Let $C=C_G(N)$, and suppose that $C\ne1$. Since
$N\cap C=Z(N)=1$, the image
\[
\overline N=NC/C\cong N
\]
is self-centralizing in $\overline G=G/C$. Indeed, if $xC$ centralizes $\overline N$, then $[x,N]\le C\cap N=1$, so $x\in C$. Any minimal normal subgroup of $\overline G$ different from $\overline N$ would centralize $\overline N$ and is therefore trivial. Hence $\overline N$ is the unique minimal normal subgroup of $\overline G$, and $\R(\overline G)=1$.

The group $\overline G$ is a proper quotient group which is a  product of two supersoluble groups, so minimality of $G$ gives
\[
\lambda(\overline G)\le1.
\]
Since $\R(\overline G)=1$ and $\Soc(\overline G)=\overline N$, it follows that
\[
\overline G/\overline N\cong G/NC
\]
is soluble. But both $N$ and $C$ fix every component of $N$, so the component action of $G$ factors through $G/NC$. The permutation group induced by $H$ on these factors is therefore soluble. This contradiction proves that $C_G(N)=1$.

If $D$ were a second minimal normal subgroup of $G$, then $[D,N]=1$,
contrary to $C_G(N)=1$. Thus $N$ is the unique minimal normal subgroup and
$E=N$. By the definition of $K$, the action of $H=M/K$ on the simple
factors of $N$ is faithful. The embedding of $K/N$ in $\Out(S)^r$ follows
from the faithful conjugation action of $K$ on $N$ modulo inner automorphisms.
\end{proof}

\section{The defining prime of the socle}\label{sec:prime}

We continue with the minimal counterexample of Section~\ref{sec:length-minimal}.
By Lemma~\ref{lem:active-lower},
\[
N=\Soc(G)=S_1\times\cdots\times S_r,
\qquad S_i\cong S=\PSL_2(p),
\]
where $p\ge5$ is prime, $C_G(N)=1$, and the action of $H=M/K$ on the
simple factors of $N$ is nonsoluble. Let $\Omega=\{S_1,\ldots,S_r\}$,
and let $\pi_i:N\longrightarrow S_i$ be the coordinate projection. Set
\[
A_N=A\cap N,
\qquad B_N=B\cap N,
\qquad F=AN\cap BN.
\]
By Lemma~\ref{lem:factorizer}, $F$ normalizes every $S_i$. Since $C_G(N)=1$,
conjugation embeds $F$ in $\prod_i\Aut(S_i)$, with $N$ mapping onto
$\prod_i\Inn(S_i)$. Consequently
\[
F/N\le\Out(S)^r\cong C_2^r.
\]
Comparing the product formulas for $G=AB$ and $G/N=(AN/N)(BN/N)$ gives
\begin{equation}\label{eq:order-identity}
|N|=|A_NB_N|\,[F/N:(A\cap B)N/N].
\end{equation}
Here $|A_NB_N|$ denotes the cardinality of the set $A_NB_N$. In particular,
\begin{equation}\label{eq:odd-part}
|N|_{2'}=|A_NB_N|_{2'}.
\end{equation}

We relate the coordinate projections of $A_N$ and $B_N$ to the local
supersoluble subgroups. Let
\[
A_0=A\cap BN,
\qquad B_0=B\cap AN.
\]
The factorizer identity gives
\[
F=A_0B_0=A_0N=B_0N,
\]
so $A_0/A_N$ and $B_0/B_N$ are $2$-groups. For each $i$, let
$\varphi_i:F\longrightarrow\Aut(S_i)$ be the conjugation map, and identify
$S_i$ with $\Inn(S_i)$. Define
\[
U_i=\varphi_i(A_0)\cap S_i,
\qquad V_i=\varphi_i(B_0)\cap S_i.
\]
The group $X_i=\varphi_i(F)$ satisfies
\[
X_i=\varphi_i(A_0)\varphi_i(B_0)
   =\varphi_i(A_0)S_i
   =\varphi_i(B_0)S_i,
\qquad |X_i:S_i|\le2.
\]
Comparing orders therefore gives
\[
|S_i|\mid 2|U_i|\,|V_i|.
\]
Moreover, $\pi_i(A_N)\unlhd U_i$ and $\pi_i(B_N)\unlhd V_i$, and the
quotients
\[
U_i/\pi_i(A_N),
\qquad V_i/\pi_i(B_N)
\]
are $2$-groups which are sections of $A_0/A_N$ and $B_0/B_N$, respectively.

By Lemma~\ref{lem:local-psl2}, exactly one of $|U_i|,|V_i|$ is divisible
by $p$.

Let $P_A\in\Syl_p(A_N)$ and $P_B\in\Syl_p(B_N)$. These are the unique
Sylow $p$-subgroups of $A_N$ and $B_N$, respectively, and are therefore
normal in $A$ and $B$. For $P\le N$, define
\[
\operatorname{supp}(P)=\{S_i\in\Omega:\pi_i(P)\ne1\}.
\]
Since the local quotients above are $2$-groups, the alternatives for
$U_i,V_i$ give
\[
\operatorname{supp}(P_B)=\Omega\setminus\operatorname{supp}(P_A).
\]
The support of $P_A$ is $A$-invariant, and its complement is $B$-invariant.
Thus both sets are $G$-invariant. Since $N$ is a minimal normal subgroup
of $G$, the action of $G$ on $\Omega$ is transitive. After interchanging
$A$ and $B$, we therefore have
\[
P_B=1,
\qquad \pi_i(P_A)\ne1\quad(1\le i\le r).
\]
Taking $p$-parts in \eqref{eq:order-identity} now gives $|P_A|=p^r$.
Writing $P_i=\pi_i(P_A)\in\Syl_p(S_i)$, we obtain
\[
P_A=P_1\times\cdots\times P_r\cong C_p^r,
\qquad P_A\lhd A.
\]

\begin{proposition}\label{prop:Fermat}
The prime $p$ satisfies $p-1=2^a$ for some $a\ge2$.
\end{proposition}

\begin{proof}
Suppose that an odd prime $\ell$ divides $(p-1)/2$, and let $\ell^e$ be
the exact $\ell$-part of $(p-1)/2$. Since
\[
A_N\le N_N(P_A)=\prod_{i=1}^rN_{S_i}(P_i),
\qquad N_{S_i}(P_i)\cong C_p:C_{(p-1)/2},
\]
the $\ell$-part of $|A_N|$ is at most $\ell^{er}$. In every coordinate,
$p\mid|U_i|$, and Lemma~\ref{lem:local-psl2}(ii) gives
$\ell\nmid|V_i|$. Since $\pi_i(B_N)\le V_i$ for every $i$, it follows
that $|B_N|_\ell=1$. Equation~\eqref{eq:odd-part} therefore yields
\[
|A_N/P_A|_\ell=\ell^{er}.
\]
As
\[
A_N/P_A\le N_N(P_A)/P_A\cong C_{(p-1)/2}^{\,r},
\]
its Sylow $\ell$-subgroup is the full coordinate subgroup
$C_{\ell^e}^{\,r}$. Its preimage in $A_N$ is
\[
W=W_1\times\cdots\times W_r,
\qquad W_i\le S_i,
\qquad W_i\cong C_p:C_{\ell^e}.
\]
The Sylow $\ell$-subgroup of the abelian group $A_N/P_A$ is
characteristic. Since $P_A$ is characteristic in $A_N$, so is $W$.
Hence $W\lhd A$.

Each complement $C_{\ell^e}$ acts faithfully on $P_i$, so
$\F(W)=P_A$. Since $A'\le\F(A)$,
\[
[W,A]\le W\cap A'
       \le W\cap\F(A)
       \le\F(W).
\]
Thus $A$ centralizes
\[
W/\F(W)\cong C_{\ell^e}^{\,r}.
\]
Conjugation by $A$ permutes the factors $W_i$ as it permutes the $S_i$.
A nontrivial permutation of these factors acts nontrivially on
$W/\F(W)$, so $A$ fixes every $S_i$. The permutation group induced by
$G=AB$ on $\Omega$ is then an image of the supersoluble group $B$ and
is soluble. This contradicts the nonsolubility of the action of $H$ on
$\Omega$.

Therefore $(p-1)/2$ has no odd prime divisor, and $p-1$ is a power of $2$.
\end{proof}

Write $P=P_A$.

\begin{proposition}\label{prop:large-fermat}
In the minimal counterexample, $p=5$.
\end{proposition}

\begin{proof}
Suppose that $p\ne5$. By Proposition~\ref{prop:Fermat}, write
$p-1=2^a$. Since $p$ is prime, $p\ge17$ and $a\ge4$. Let
\[
m=a-1,
\qquad \nu=\frac{p+1}{2}.
\]
Then $m\ge3$ and $\nu$ is odd. Lemma~\ref{lem:local-psl2}(ii)--(iii) gives
\[
U_i\le N_{S_i}(P_i)\cong C_p:C_{2^m},
\qquad V_i\le N_{S_i}(T_i),
\]
where $T_i$ is a nonsplit torus of $S_i$ of order $\nu$. Consequently
\[
A_N/P\le C_{2^m}^{\,r},
\qquad B_N\le\prod_{i=1}^rN_{S_i}(T_i).
\]
Let
\[
T=T_1\times\cdots\times T_r\cong C_\nu^{\,r}.
\]
The odd part of $|A_N|$ is $p^r$, whereas $p\nmid|B_N|$. Thus
\eqref{eq:odd-part} implies $|B_N|_{2'}=\nu^r$. As $N_{S_i}(T_i)/T_i\cong C_2$ for every $i$, the index $|B_N:B_N\cap T|$ is a power of $2$. Hence
\[
|B_N\cap T|=\nu^r=|T|,
\]
so $T\le B_N$ and $B_N/T\le C_2^r$. Write
\[
|A_N/P|=2^x,
\qquad |B_N/T|=2^y.
\]

We have
\[
\F(A_N)=P,
\qquad \F(B_N)=T.
\]
Indeed,
\[
C_{\prod_iN_{S_i}(P_i)}(P)=P.
\]
The normal subgroups $O_2(A_N)$ and $P$ centralize one another, so this
equality forces $O_2(A_N)=1$. Since $A_N/P$ is a $2$-group, it follows
that $\F(A_N)=P$. Similarly,
\[
C_{\prod_iN_{S_i}(T_i)}(T)=T.
\]
The normal subgroups $O_2(B_N)$ and $T$ centralize one another, giving
$O_2(B_N)=1$ and hence $\F(B_N)=T$.

Since $A_N\lhd A$, $B_N\lhd B$, and the derived subgroup of a
supersoluble group is nilpotent,
\[
[A_N,A]\le A_N\cap A'
         \le A_N\cap\F(A)
         \le\F(A_N)=P,
\]
and similarly $[B_N,B]\le T$. Thus $A$ centralizes $A_N/P$ and $B$
centralizes $B_N/T$.

Let $o_A$ and $o_B$ be the numbers of $A$-orbits and $B$-orbits on
$\Omega$. A tuple in $A_N/P\le C_{2^m}^{\,r}$ fixed by $A$ is determined
by its value at one coordinate in each $A$-orbit: conjugation carries
each coordinate factor isomorphically to the corresponding factor.
The same observation applies to $B_N/T\le C_2^r$. Therefore
\[
x\le mo_A,
\qquad y\le o_B.
\]
Taking $2$-parts in \eqref{eq:order-identity} and using $|N|_2=2^{ar}$
and $|F/N|\le2^r$ gives
\[
2^{ar}
\le |A_NB_N|_2\,2^r
\le 2^{x+y+r}.
\]
Hence $x+y\ge(a-1)r=mr$, and thus
\begin{equation}\label{eq:length-orbit-1}
mr\le mo_A+o_B.
\end{equation}

Form a bipartite multigraph with the $A$-orbits and the $B$-orbits on
$\Omega$ as its two vertex sets. Each $S_i$ gives an edge joining the
two orbits containing it. The action of $G=\langle A,B\rangle$ on
$\Omega$ is transitive, so this graph is connected. It has $r$ edges
and $o_A+o_B$ vertices, whence
\begin{equation}\label{eq:length-orbit-2}
o_A+o_B\le r+1.
\end{equation}
Combining \eqref{eq:length-orbit-1} and \eqref{eq:length-orbit-2} gives
\[
mr\le(m-1)o_A+r+1,
\qquad (m-1)(r-o_A)\le1.
\]
Since $m\ge3$, we obtain $o_A=r$. Thus $A$ fixes every component of
$N$, and the permutation group induced by $G=AB$ on $\Omega$ is a
soluble image of $B$. This again contradicts the nonsolubility of the
action of $H$ on $\Omega$.
\end{proof}

\section{Cycle bounds}\label{sec:cycles}

The case $p=5$ requires bounds on the numbers of cycles of elements of the upper simple section. For a permutation $x$ of a finite set $\Omega$, let $o_\Omega(x)$ denote its number of cycles.

\begin{lemma}\label{lem:cycle-prime}
Let $H$ act faithfully on a finite set $\Omega$, and suppose that its soluble radical $R=\R(H)$ is transitive and $H/R$ is nonabelian simple.
\begin{enumerate}[label=\textup{(\roman*)}]
\item If $s$ is prime and $x\in H$ is an $s$-element with $xR\ne R$, then
\[
o_\Omega(x)\le \frac{2s-1}{s^2}|\Omega|.
\]
\item If $x\in H$ is a $2$-element and $xR$ has order divisible by $4$, then
\[
o_\Omega(x)\le\frac58|\Omega|.
\]
\end{enumerate}
\end{lemma}

\begin{proof}
We prove both assertions by induction on $|\Omega|$. In \textup{(i)} let $\beta=(2s-1)/s^2$, and in \textup{(ii)} let $s=2$ and $\beta=5/8$. In both cases $1/s\le\beta$. Choose a minimal normal subgroup $U$ of $H$ contained in $R$. Then $U$ is elementary abelian.

Suppose first that $U$ is transitive. Since $U$ is abelian and the action is faithful, $U$ is regular. If $x$ has no fixed points, then every cycle has length at least $s$, so $o_\Omega(x)\le|\Omega|/s$. Otherwise we may conjugate $x$ by an element of $U$ so that it fixes the identity point of $U$. Since a regular abelian subgroup is self-centralizing in the full symmetric group, $C_H(U)=U$. Thus $x$ induces a nontrivial linear transformation on $U$, and
\[
|U:C_U(x)|\ge s.
\]
In characteristic $s$, this follows from the Jordan form of a nontrivial unipotent transformation. In characteristic $t\ne s$, the direct sum of the nontrivial irreducible constituents has dimension $e$, and some irreducible constituent of dimension $j\le e$ has an eigenvalue whose order is divisible by $s$. Hence $s\mid t^j-1$, and $t^e\ge t^j\ge s+1$.

Under the hypothesis of \textup{(ii)}, the linear transformation induced by $x$ has order divisible by $4$. If $U$ has characteristic $2$, some Jordan block has size at least $3$, so the fixed-space codimension is at least $2$. In odd characteristic $t$, choose an eigenvalue of order divisible by $4$, and let $j$ be its degree over $\mathbb F_t$. Then $4\mid t^j-1$, so $|U:C_U(x)|\ge t^j\ge5$. In either case,
\[
|U:C_U(x)|\ge4.
\]
The fixed points of $x$ are precisely $C_U(x)$, and every other cycle has length at least $s$. These bounds give
\[
o_\Omega(x)
\le |C_U(x)|+\frac{|U|-|C_U(x)|}{s}
\le\beta|U|.
\]

Now suppose that $U$ is intransitive, and let $\mathcal B$ be the set of its orbits. Let $C$ be the kernel of the action of $H$ on $\mathcal B$. Since $H/R$ is simple, either $C\le R$ or $CR=H$. If $C\le R$, then $\R(H/C)=R/C$ is transitive on $\mathcal B$, and the quotient by this radical is again $H/R$. The image of $x$ has the same order modulo the radical as $xR$. Induction gives at most $\beta|\mathcal B|$ cycles on $\mathcal B$. A cycle through blocks of size $b$ contributes at most $b$ cycles on $\Omega$, proving the required bound.

Assume that $CR=H$. For a block $\Delta\in\mathcal B$ fixed by $x$, let $H_\Delta$ and $R_\Delta$ denote its setwise stabilizers in $H$ and $R$. Then
\[
H_\Delta=CR_\Delta,
\qquad
H_\Delta/R_\Delta\cong C/(C\cap R)\cong H/R.
\]
The subgroup $R_\Delta$ is soluble and normal in $H_\Delta$, and it is transitive on $\Delta$ because it contains $U$. Let $J$ be the kernel of the action of $H_\Delta$ on $\Delta$.

The image of $C$ on $\Delta$ is nonsoluble. Indeed, the images of $C$ on the different $U$-orbits are conjugate under $H$, since $C\lhd H$. Their kernels have trivial intersection because $H$ acts faithfully on $\Omega$. If one image were soluble, all would be soluble, and the embedding of $C$ in their direct product would make $C$ soluble, contrary to $C/(C\cap R)\cong H/R$. Hence $H_\Delta/J$ is nonsoluble. If $JR_\Delta=H_\Delta$, this quotient would be an image of the soluble group $R_\Delta$. Simplicity of $H_\Delta/R_\Delta$ therefore gives $J\le R_\Delta$.

Thus the faithful induced group on $\Delta$ has soluble transitive radical $R_\Delta/J$ and quotient isomorphic to $H/R$. The natural isomorphism $H_\Delta/R_\Delta\cong H/R$ sends $xR_\Delta$ to $xR$, so the hypothesis on the order of this image is unchanged. Induction applies on every fixed block.

Suppose that the $n$ blocks have size $b$ and that $x$ fixes $f$ of them. The fixed blocks contribute at most $\beta bf$ cycles. Every nonfixed block orbit has length at least $s$ and contributes at most $b$ cycles. Hence
\[
o_\Omega(x)\le \beta bf+\frac{b(n-f)}s
\le\beta bn=\beta|\Omega|.\qedhere
\]
\end{proof}

\begin{lemma}\label{lem:upper-elements}
Let
\[
T=\PSL_2(q)\le X\le\PGL_2(q),
\]
where $q\ge5$ is prime. Suppose that
\[
X=UV
\]
with $U,V\le X$ supersoluble. After interchanging $U$ and $V$, the subgroup $U\cap T$ contains an element of order $q$. Moreover, one of the following holds:
\begin{enumerate}[label=\textup{(\roman*)}]
\item $V\cap T$ contains an element of odd prime order $s$ dividing $q+1$;
\item $q+1$ is a power of $2$ and $V\cap T$ contains an element of order $4$.
\end{enumerate}
\end{lemma}

\begin{proof}
Let $U_0=U\cap T$, $V_0=V\cap T$, and
\[
d=|X:T|,\qquad a=|U:U_0|,\qquad b=|V:V_0|.
\]
Since $d\le2$ and the images of $U$ and $V$ multiply to $X/T$, we have $ab/d\in\{1,2\}$. The product formula gives
\[
|T|\,|U\cap V|=\frac{ab}{d}|U_0||V_0|,
\]
and hence
\[
|T|\mid2|U_0||V_0|.
\]
After interchanging the factors, we may assume that $q\mid |U_0|$. Since $q$ is the largest prime divisor of $|T|$, the Sylow $q$-subgroup $P$ of the supersoluble group $U_0$ is normal. Thus
\[
U_0\le N_T(P)\cong C_q:C_{(q-1)/2}.
\]
If an odd prime $s$ divides $q+1$, then $(s,2|N_T(P)|)=1$, so the preceding divisibility forces $s\mid |V_0|$. This proves \textup{(i)}.

Suppose that $q+1=2^e$. Then $(q-1)/2=2^{e-1}-1$ is odd, and
\[
|V_0|_2\ge2^{e-1}.
\]
A Sylow $2$-subgroup of $T$ is dihedral of order $2^e$ by \cite[Chap.~II, Theorem 8.27]{Huppert}. If $q\ne7$, then $e\ge5$, and every subgroup of this dihedral group of order at least $2^{e-1}$ contains an element of order $4$.

Let $q=7$. The subgroup $U$ normalizes its Sylow $7$-subgroup and hence lies in the point stabilizer $N_X(P)$ in the natural action of $X$ of degree $8$. Since $X=UV$, the group $V$ is transitive on these eight points. If $X=T$, then $8\mid |V|$, so $V$ contains a Sylow $2$-subgroup of $T$, which is dihedral of order $8$.

Suppose that $X=\PGL_2(7)$. If $7\mid |V|$, then normality of its Sylow $7$-subgroup places $V$ in a point stabilizer, contrary to transitivity. If $3\mid |V|$, then a Sylow $3$-subgroup $C$ of $V$ is normal, since $7\nmid |V|$. Thus $V\le N_X(C)$, whose order is $12$, again contradicting transitivity of degree $8$. Hence $V$ is a $2$-group and lies in a Sylow $2$-subgroup $P_2\cong D_{16}$ of $X$, with $|V|\ge8$. Every subgroup of index $2$ in $P_2$ contains the Frattini subgroup $\Phi(P_2)\cong C_4$. This applies to both $P_2\cap T$ and $V$ when $|V|=8$, and the case $V=P_2$ is immediate. Therefore $\Phi(P_2)\le V\cap T$.
\end{proof}

\begin{lemma}\label{lem:simple-block-cycles}
Let $S=\PSL_2(q)$, where $q\ge5$ is prime, and let $S$ act faithfully and transitively on a finite set $\mathcal B$. For $x,y\in S$, the following bounds hold.
\begin{enumerate}[label=\textup{(\roman*)}]
\item If $x$ has order $q$, then
\[
o_{\mathcal B}(x)\le \frac{2}{q+1}|\mathcal B|.
\]
\item If $s$ is an odd prime dividing $q+1$ and $y$ has order $s$, then
\[
o_{\mathcal B}(y)\le \frac35|\mathcal B|.
\]
\item If $q+1$ is a power of $2$ and $y$ has order $4$, then
\[
o_{\mathcal B}(y)\le \frac47|\mathcal B|.
\]
\end{enumerate}
\end{lemma}

\begin{proof}
For $z\in S$, let $\operatorname{fpr}_{\mathcal B}(z)$ be the proportion of points of $\mathcal B$ fixed by $z$. If $J$ is a point stabilizer, then
\[
\operatorname{fpr}_{\mathcal B}(z)
 =\frac{|z^S\cap J|}{|z^S|}.
\]
We use Dickson's subgroup classification; see \cite[Chap.~II, Theorem  8.27]{Huppert}.

Let $x$ have order $q$. A proper subgroup containing $x$ lies in a Borel subgroup. The class of $x$ has size $(q^2-1)/2$, while a Borel contains $(q-1)/2$ elements from this class. Hence
\[
\operatorname{fpr}_{\mathcal B}(x)\le\frac1{q+1}.
\]
Every nonfixed cycle has length $q$, so
\[
\frac{o_{\mathcal B}(x)}{|\mathcal B|}
\le \frac1q+\left(1-\frac1q\right)\frac1{q+1}
=\frac2{q+1}.
\]

Now let $y$ have odd prime order $s$ dividing $q+1$. Then $y$ is nonsplit semisimple and $|y^S|=q(q-1)$. A Borel subgroup or a split-torus normalizer contains no element of order $s$, while a nonsplit-torus normalizer contains at most two elements of $y^S$. The remaining proper subgroup types are contained in $A_4$, $S_4$, or $A_5$. If $q=5$, then $s=3$, and every proper subgroup of $A_5$ contains at most eight elements of order $3$, giving the bound $8/20=2/5$. If $q=7$, no odd prime divides $q+1$. For $q\ge11$, each of $A_4,S_4,A_5$ contains at most $24$ elements of any fixed odd prime order, and
\[
\frac{24}{q(q-1)}\le\frac{24}{110}<\frac25.
\]
Thus $\operatorname{fpr}_{\mathcal B}(y)\le2/5$. Since every nonfixed cycle has length $s$,
\[
\frac{o_{\mathcal B}(y)}{|\mathcal B|}
\le \frac1s+\left(1-\frac1s\right)\frac25
\le\frac35.
\]

Finally suppose that $q+1$ is a power of $2$ and $|y|=4$. Then $q\equiv7\pmod8$. A proper subgroup containing $y$ lies in a nonsplit-torus normalizer or in an $S_4$. The class of $y$ has size $q(q-1)$; the torus normalizer meets it in at most two elements and $S_4$ contains six elements of order $4$. Hence
\[
\operatorname{fpr}_{\mathcal B}(y)\le\frac6{q(q-1)}\le\frac17.
\]
Every nonfixed cycle has length at least $2$, so
\[
\frac{o_{\mathcal B}(y)}{|\mathcal B|}
\le \frac12+\frac12\operatorname{fpr}_{\mathcal B}(y)
\le\frac47.\qedhere
\]
\end{proof}

\section{Proof of the nonsoluble-length theorem}\label{sec:length-a5}

\begin{proof}[Proof of Theorem~\ref{thm:main}]
Suppose that the theorem is false, and choose the minimal counterexample of Section~\ref{sec:length-minimal}. Keep the notation in  Section~\ref{sec:prime}. By Proposition~\ref{prop:large-fermat}, $p=5$, so $N=A_5^r$. Write
\[
P=P_A=P_1\times\cdots\times P_r\cong C_5^r.
\]
Thus $P\le A_N$ and $P\lhd A$.

For every $i$, Lemma~\ref{lem:local-psl2}(ii) gives
\[
U_i\le N_{S_i}(P_i)\cong D_{10}.
\]
Since $|S_i|\mid2|U_i||V_i|$ and $3\nmid|U_i|$, we have $3\mid|V_i|$. Also $5\nmid|V_i|$, so the Sylow $3$-subgroup $Q_i$ of $V_i$ is normal and has order $3$. Hence
\[
V_i\le N_{S_i}(Q_i)\cong S_3.
\]
Consequently,
\[
\begin{aligned}
A_N&\le\prod_{i=1}^rN_{S_i}(P_i)\cong D_{10}^r,\\
B_N&\le\prod_{i=1}^rN_{S_i}(Q_i)\cong S_3^r.
\end{aligned}
\]
Let $Q=Q_1\times\cdots\times Q_r\cong C_3^r$. Since $3\nmid|A_N|$, equation~\eqref{eq:odd-part} gives $|B_N|_3=3^r$. The quotient $\bigl(\prod_iN_{S_i}(Q_i)\bigr)/Q$ is a $2$-group, so
\[
|B_N\cap Q|=3^r=|Q|.
\]
Thus $Q\le B_N$. It is the unique Sylow $3$-subgroup of $B_N$, and therefore characteristic in $B_N$ and normal in $B$. Let
\[
D_A=A_N/P\le C_2^r,\qquad
D_B=B_N/Q\le C_2^r,
\]
and write $c=\dim_{\mathbb F_2}D_A$ and $d=\dim_{\mathbb F_2}D_B$. Taking $2$-parts in~\eqref{eq:order-identity} and using $F/N\le C_2^r$, we obtain
\begin{equation}\label{eq:length-cycle-1}
c+d\ge r.
\end{equation}

The centralizer equalities
\[
C_{D_{10}^r}(P)=P,\qquad C_{S_3^r}(Q)=Q
\]
give $\F(A_N)=P$ and $\F(B_N)=Q$. Indeed, any normal $2$-subgroup of either group centralizes its normal Hall subgroup of odd order and is therefore trivial. Since a supersoluble group has nilpotent derived subgroup,
\[
[A_N,A]\le A_N\cap A'
   \le A_N\cap\F(A)
   \le\F(A_N)=P,
\]
and similarly $[B_N,B]\le Q$. Hence $A$ centralizes $D_A$ and $B$ centralizes $D_B$.

For $x\in A$, an element of $D_A\le C_2^r$ fixed by $x$ is determined by one coordinate on each cycle of $x$ on $\Omega$. Thus
\begin{equation}\label{eq:length-cycle-2}
c\le o_\Omega(x)\qquad(x\in A).
\end{equation}
Similarly,
\begin{equation}\label{eq:length-cycle-3}
d\le o_\Omega(y)\qquad(y\in B).
\end{equation}

Let $\overline G=G/K$. The group $H=M/K$ is normal in $\overline G$ and acts faithfully on $\Omega$, with soluble radical $R=L/K$ and
\[
H/R\cong\PSL_2(q).
\]
Since $\overline G$ is transitive on $\Omega$, it permutes the $H$-orbits transitively. For an $H$-orbit $\Delta$, let $J_\Delta$ be the kernel of the action of $H$ on $\Delta$. The induced actions are conjugate under $\overline G$, so the groups $H/J_\Delta$ are isomorphic. Moreover, $H$ embeds in their direct product. Since $H$ is nonsoluble, every $H/J_\Delta$ is nonsoluble.

Now $J_\Delta R/R\lhd H/R$. If $J_\Delta R=H$, then $H/J_\Delta$ would be an image of the soluble group $R$, a contradiction. Thus $J_\Delta\le R$. It follows that the faithful induced group on $\Delta$ has soluble radical $R/J_\Delta$ and quotient
\[
(H/J_\Delta)/(R/J_\Delta)\cong H/R.
\]

If $R$ is transitive on $\Delta$, Lemma~\ref{lem:cycle-prime} applies to this induced group. If $R$ is intransitive, let $\mathcal B$ be the set of its orbits on $\Delta$, and let $C$ be the kernel of the action of $H$ on $\mathcal B$. Then $R\le C$, and $C/R\lhd H/R$. Since $H$ is transitive on $\mathcal B$ and $|\mathcal B|>1$, we have $C\ne H$, so simplicity gives $C=R$. Thus $H/R$ acts faithfully and transitively on $\mathcal B$. If the blocks have size $b$, each cycle on $\mathcal B$ contributes at most $b$ cycles on $\Delta$. The bounds in Lemma~\ref{lem:simple-block-cycles} therefore also bound the corresponding cycle proportions on $\Delta$.

By Proposition~\ref{prop:upper-factorized}, Lemma~\ref{lem:upper-elements} applies to the images of $A$ and $B$ in $G/L$, with $T=M/L$. Its elements lie in the images of $A\cap M$ and $B\cap M$, since
\[
AL\cap M=(A\cap M)L,\qquad
BL\cap M=(B\cap M)L.
\]
For a prime $\ell$, every $\ell$-element of either quotient intersection has an $\ell$-element lift in the corresponding intersection, because a Sylow $\ell$-subgroup of the preimage of its cyclic subgroup maps onto it.

We may therefore choose $x,y\in M$, with $x\in A\cap M$ and $y\in B\cap M$ or with the roles of $A$ and $B$ reversed, such that $x$ is a $q$-element and $xL$ has order $q$. In case~\textup{(i)} of Lemma~\ref{lem:upper-elements}, the element $y$ is an $s$-element and $yL$ has order $s$, for an odd prime $s\mid q+1$. In case~\textup{(ii)}, the element $y$ is a $2$-element and $yL$ has order $4$. Write
\[
\widehat x=xK,\qquad \widehat y=yK
\]
for their images in $H$. Under $H/R\cong M/L$, the cosets $\widehat xR$ and $\widehat yR$ correspond to $xL$ and $yL$. Moreover, $x$ and $\widehat x$ induce the same permutation of $\Omega$, as do $y$ and $\widehat y$.

Suppose first that case~\textup{(i)} holds. On every $H$-orbit $\Delta$,
\[
o_\Delta(\widehat x)\le\frac{2q-1}{q^2}|\Delta|,
\qquad
o_\Delta(\widehat y)\le\frac35|\Delta|.
\]
For $\widehat x$, use Lemma~\ref{lem:cycle-prime}(i) when $R$ is transitive and Lemma~\ref{lem:simple-block-cycles}(i) otherwise, together with
\[
\frac2{q+1}\le\frac{2q-1}{q^2}.
\]
For $\widehat y$, the transitive case gives $(2s-1)/s^2\le5/9<3/5$, and the intransitive case follows from Lemma~\ref{lem:simple-block-cycles}(ii). Summing over the $H$-orbits and using \eqref{eq:length-cycle-2}--\eqref{eq:length-cycle-3}, we obtain
\[
c+d\le o_\Omega(x)+o_\Omega(y)
\le\left(\frac{2q-1}{q^2}+\frac35\right)r
\le\frac{24}{25}r<r,
\]
contrary to \eqref{eq:length-cycle-1}.

Suppose that case~\textup{(ii)} holds. Then $q+1$ is a power of $2$ and $q\ge7$. The same bound applies to $\widehat x$. For $\widehat y$, Lemma~\ref{lem:cycle-prime}(ii) gives the bound $5/8$ on an $H$-orbit on which $R$ is transitive, while Lemma~\ref{lem:simple-block-cycles}(iii) gives the stronger bound $4/7$ when $R$ is intransitive. Hence
\[
c+d\le o_\Omega(x)+o_\Omega(y)
\le\left(\frac{2q-1}{q^2}+\frac58\right)r
\le\frac{349}{392}r<r,
\]
again contradicting \eqref{eq:length-cycle-1}. Therefore no minimal counterexample exists, and $\lambda(G)\le1$.
\end{proof}

\section*{Acknowledgements}
Li acknowledges support from the NSFC  (No. 12471017) and the Scientific Research Foundation for Advanced Talents of Suqian University (No. 2022XRC069).  Yang acknowledges support from the Simons Foundation (Award No.~918096).   During  this work,  the authors used OpenAI's ChatGPT  for mathematical exploration on the nonsoluble  length of the product of two supersoluble  groups  and for  editorial assistance. The authors reviewed and edited the content as required  and take full responsibility for the content.

\section*{Disclosure Statement}
The authors declare no competing interests.

\section*{Data Availability Statement}
Data sharing is not applicable to this article, as no data sets were generated or analysed during the current study.

\end{document}